\documentclass[11pt]{article}
\usepackage{amsmath}
\usepackage{amsfonts}  
\usepackage{amssymb}
\usepackage[utf8]{inputenc}         
\usepackage[T1]{fontenc}              
\usepackage{color}
\usepackage{tikz}
\usetikzlibrary{arrows, decorations.markings, decorations.pathmorphing, backgrounds, positioning, fit, petri}
\newcommand{\qed}{\hspace*{\fill}$\Box $\\ \vspace{0.5cm}} 
\newcommand{\bewanf}{\noindent {\bf Proof:}\quad }    
\newcommand{\satzanf}{\begin{samepage}\begin{satz}}
\newcommand{\satzende}{\end{satz}\end{samepage}}
\newcommand{\koranf}{\begin{samepage}\begin{kor}}
\newcommand{\korende}{\end{kor}\end{samepage}}
\newcommand{\lemanf}{\begin{samepage}\begin{lem}}
\newcommand{\lemende}{\end{lem}\end{samepage}}
\newcommand{\bspanf}{\begin{beispiel}}
\newcommand{\bspende}{\end{beispiel}}
\newcommand{\defanf}{\begin{definition}}
\newcommand{\defende}{\end{definition}}
\newcommand{\propanf}{\begin{proposition}}
\newcommand{\propende}{\end{proposition}}
\newtheorem{satz}{Theorem}[section]             
\newtheorem{lem}[satz]{Lemma}
\newtheorem{beispiel}[satz]{Example}
\newtheorem{definition}[satz]{Definition}
\newtheorem{proposition}[satz]{Proposition}
\newtheorem{kor}[satz]{Corollary}
\newcommand{\be}{\begin{equation}}
\newcommand{\ee}{\end{equation}}
\newcommand{\bl}[1]{\begin{equation}\label{#1}}
\newcommand{\ben}{\begin{enumerate}}
\newcommand{\een}{\end{enumerate}}
\begin{document}
\title{The Local Structure of Injective LOT-Complexes}
\author{Jens Harlander and Stephan Rosebrock}
\maketitle
\thispagestyle{empty}
\begin{abstract}
Labeled oriented trees, LOT's, encode spines of ribbon discs in the 4-ball and ribbon 2-knots in the 4-sphere. The unresolved asphericity question for these spines is a major test case for Whitehead's asphericity conjecture. In this paper we give a complete description of the link of a reduced injective LOT complex. An important case is the following: If $\Gamma$ is a reduced injective LOT that does not contain boundary reducible sub-LOTs, then $lk(K(\Gamma))$ is a bi-forest. As a consequence $K(\Gamma)$ is aspherical, in fact DR, and its fundamental group is locally indicable. We also show that a general injective LOT complex is aspherical. Some of our results have already appeared in print over the last two decades and are collected here.
\end{abstract}
{\sc Keywords:} Labeled oriented tree; non-positive immersion; coloring test; weight test; locally indicable\\

{\sc MSC 2020:} 57K20; 57M07; 20F06; 20F65; 20F67\\
\section{Introduction}

A labeled oriented graph (LOG) $\Gamma = (E, V, s, t, \lambda)$ consists of two sets $E$, $V$ of edges and vertices, and three maps $s, t, \lambda\colon E\to V$ called, respectively source, target and label. $\Gamma$ is said to be a labeled oriented forest (LOF) when the underlying graph is a forest, and it is called labeled oriented tree (LOT) if the underlying graph is a tree. The associated LOG presentation is defined as
$$P(\Gamma)=\langle V\ |\  s(e)\lambda(e)=\lambda(e)t(e),\ e\in E \rangle.$$
The LOG complex $K(\Gamma)$ is the standard 2-complex defined by the presentation, and $G(\Gamma)=\pi_1(K(\Gamma))$. 

It is known that LOT-complexes are spines of ribbon 2-knot complements \cite{How83}. So the study of LOTs and LOFs is an extension of classical knot and link theory. Asphericity, known for classical knots, is unresolved for LOTs. The asphericity question for LOTs is of central importance to Whitehead's asphericity conjecture: The subcomplex of an aspherical 2-complex is aspherical. See Berrick/Hillman \cite{BerrickHillman}, Bogley \cite{Bogley}, and Rosebrock \cite{Ro18}.

A sub-LOG $\Gamma_0=(E_0, V_0)\subseteq\Gamma$ is a subgraph so that $E_0\ne\emptyset$ and $\lambda\colon E_0\to V_0$. A LOG is called {\em boundary reduced} if whenever $v$ is a vertex of valency 1 then $v=\lambda(e)$ for some edge $e$. It is called {\em interior reduced} if for every vertex $v$ no two edges starting or terminating at $v$ carry the same label. It is called {\em compressed }if for every edge $e$ the label $\lambda(e)$ is not equal to $s(e)$ or $t(e)$. Finally, a LOG is {\em reduced} if it is boundary reduced, interior reduced, and compressed. Given a LOG, reductions can be performed to produce a reduced LOG, and, in case the LOG is a LOF, this process does not affect the homotopy type of the LOF complex. A LOG is called {\em injective} if the labeling map $\lambda\colon E\to V$ is injective.

The main goal of this paper is to give an explicit description of the link $lk(K(\Gamma))$, the boundary of a regular neighborhood of the unique vertex, in case $\Gamma$ is a reduced injective LOF. 
This description is of fundamental importance. A graph $\Lambda$ is called a {\em bi-forest} if it consists of disjoint forests $A$ and $B$ so that edges that are contained in neither of $A,B$ connect a vertex from $A$ to a vertex from $B$. Assume $\Gamma$ is a reduced injective LOF without boundary reducible sub-LOTs. We will show:
\begin{itemize}
    \item $lk(K(\Gamma))$ is a bi-forest (Theorem \ref{thm:main}).
    \item $K(\Gamma)$ admits a zero/one-angle structure that satisfies the coloring test. Therefore $K(\Gamma)$ is aspherical, has collapsible non-positive immersions, and $G(\Gamma)$ is locally indicable (Theorem \ref{thm:main} together with Theorem  \ref{slfwt}). 
\end{itemize}

We will explain the vocabulary in detail in the next section. We also address the situation for general reduced injective LOFs. Assume $\Gamma$ is a reduced injective LOF with disjoint sub-LOTs $\Gamma_1,\ldots, \Gamma_m$, and assume that the quotient LOF $\bar\Gamma$ is without boundary reducible sub-LOTs. Then
\begin{itemize}
    \item $lk(K(\Gamma))$ is a relative bi-forest (Theorem \ref{thm:rellbf}).
    \item $K(\Gamma)$ admits a zero/one-angle structure that satisfies the relative coloring test (Theorem \ref{thm:relcolLOT}).
\end{itemize}
This can be used to show
\begin{itemize}
    \item If $\Gamma$ is an injective LOF then $K(\Gamma)$ is aspherical (Theorem \ref{thm:aspherical}).
\end{itemize}

We believe it is true that $K(\Gamma)$ has collapsible non-positive immersions when $\Gamma$ is reduced and injective. But this is work in progress.\\

Some history. Reorientations of LOFs in the study of the asphericity question goes back to Huck/Rosebrock \cite{HR95}. The fact that reduced injective LOFs without boundary reducible sub-LOTs satisfy the coloring test was shown in \cite{HR01} by the same authors. Relative versions of the Huck/Rosebrock techniques were developed by the authors in \cite{HR17} and \cite{HR21}. As a consequence it was shown that injective LOT complexes are VA and therefore aspherical. Many of the results and techniques presented in this paper are indeed not new. The explicit (relative) bi-forest description of the link of a LOT complex $K(\Gamma)$, where $\Gamma$ is reduced and injective, is new. The results concerning non-positive immersions are new and are a consequence of work of Wise. The purpose of this paper is to: (1) collect fundamentally important techniques and results in a single place; (2) present new and known theorems in a unified and transparent way which allows for streamlined and simplified proofs of known results; (3) provide context that fits with the narrative of geometric group theory as it emerged over the last 20 years.

\section{Coloring tests}
A map between 2-complexes is called {\em combinatorial} if it maps open cells homeomorphically to open cells. A 2-complex is called {\em combinatorial} if the attaching maps for the 2-cells are combinatorial. Throughout the paper all 2-complexes are combinatorial.

Let $K$ be a 2-complex. If we assign numbers $\omega(c)$ to the corners $c$ of the 2-cells of $K$ we arrive at an {\em angled 2-complex} (see Gersten \cite{Ger87}). Curvature in an angled 2-complex is defined in the following way. If $v$ is a vertex of $K$ then $\kappa(v)$, the curvature at $v$, is
$$\kappa(v)=2-\chi(lk(v))-\sum \omega(c_i),$$ where the sum is taken over all the corners at $v$. If $d$ is a 2-cell of $K$ then $\kappa(d)$, the curvature of $d$, is 
$$\kappa(d)=\sum \omega(c_j)-(|\partial d|-2),$$ where the sum is taken over all the corners in $d$ and $|\partial d|$ is the number of edges in the boundary of the 2-cell. The combinatorial Gauss-Bonnet Theorem states that
$$2\chi(K)=\sum_{v\in K} \kappa(v) + \sum_{d\in K} \kappa(d).$$
This was first proven by Ballmann and Buyalo \cite{BallBuy}, and later observed by McCammond and Wise \cite{McCammondWise}.
Let $X\to K$ be a combinatorial map between 2-complexes. Note that if $K$ is an angled 2-complex then the angles in the 2-cells of $K$ can be pulled back to make $X$ into an angled 2-complex. We call this the angle structure on $X$ {\em induced} by the combinatorial map. 

A map $X\to K$ is an {\em immersion} if it is a local injection. 

\defanf A 2-complex $K$ has {\em collapsing non-positive immersions} if for every combinatorial immersion $X\to K$, where $X$ is finite, connected, either $\chi(X)\le 0$ or $X$ collapses to a point.
\defende

The concept of non-positive immersions is due to Wise. An unpublished preprint was available on his website since 1996. See \cite{Wise22} for a recent publication on the topic. Among other things Wise showed that if $K$ has non-positive immersions then $\pi_1(K)$ is locally indicable.

An angled 2-complex where all angles are either 0 or 1 is called a {\em zero/one-angled 2-complex}.
The following coloring test is due to Sieradski \cite{S83}. See also Gersten \cite{Ger87} for background and a more general weight test.

\defanf(Coloring test)
Let $K$ be a zero/one-angled 2-complex. Then $K$ {\em satisfies the coloring test} if 
\begin{enumerate}
    \item the curvature of every 2-cell is $\le 0$;
    \item for every vertex $v$: If $c_1\cdots c_n$ is a simple reduced cycle in $lk(v)$, then\\ $2-\sum_{i=1}^n \omega(c_i)\le 0$. 
\end{enumerate}
\defende

\defanf
Let $\Lambda$ be a graph. 
\begin{itemize}
    \item A cycle of edges $e_1\cdots e_n$ in $\Lambda$ is {\em reduced} if there does not exist an $e_i$ so that $e_{i+1}=\bar e_i$, where $\bar e_i$ is the edge $e_i$ with reversed orientation. A cycle of edges $e_1\cdots e_n$ in $\Lambda$ is {\em homology reduced} if there does not exist a pair $e_i, e_j$ so that $e_j=\bar e_i$.
\end{itemize}    
Let $K$ be a 2-complex. 
\begin{itemize}
    \item A {\em spherical diagram} is a combinatorial map $f\colon S\to K$ where $S$ is a cell structure on the 2-sphere $S^2$. It is {\em reduced} if, for every vertex $v\in S$, $f$ maps $lk(v)$ to a reduced cycle. It is {\em vertex reduced} if, for every vertex $v\in S$, $f$ maps $lk(v)$ to a homology reduced cycle. 
    \item $K$ is {\em diagrammatically reducible (DR)} if there do not exist reduced spherical diagrams over $K$. $K$ is {\em vertex aspherical (VA)} if there do not exist vertex reduced spherical diagrams over $K$. 
    \end{itemize}

Let $(L,K)$ be a 2-complex pair. 
\begin{itemize}
    \item The pair is {\em relatively DR} if every reduced spherical diagram over $L$ is a diagram over $K$. It is {\em relatively VA} if every vertex reduced spherical diagram over $L$ is a diagram over $K$.
\end{itemize}
\defende

\satzanf\label{thm:wise}
Let $K$ be a 2-complex. If $K$ satisfies Gersten’s weight test, it
is DR and therefore aspherical. In particular this holds when it is a zero/one angled
2-complex which satisfies the coloring test. Moreover, in that case it also
has collapsing non-positive immersions, which implies that $\pi_1(K)$ is locally indicable.\satzende

The DR part was shown by Sieradski \cite{S83} (for zero/one-angled complexes) and generalized
by Gersten \cite{Ger87} (for angled 2-complexes).  The non-positive immersions part
was proved by Wise \cite{Wise04} (In fact, Wise proved that it has nonpositive sectional curvature, which implies collapsing non-positive immersions). 

\defanf
Let $\Lambda$ be a graph and $\Lambda'$ be a subgraph, both are allowed to be disconnected. We say {\em $\Lambda$ is a forest relative to $\Lambda'$} if every homology reduced cycle in $\Lambda$ is contained in $\Lambda'$.
\defende

\defanf (Relative coloring test)
Let $K$ be a subcomplex of a zero/one-angled 2-complex $L$. Then $(L,K)$ satisfies the relative coloring test if 
\begin{enumerate}
\item the curvature of every 2-cell $d\in L-K$ is $\le 0$;
\item for every vertex $v$: If $c_1\cdots c_n$ is a homology reduced cycle in $lk(v,L)$ not entirely contained in $lk(v,K)$, then $$2-\sum_{i=1}^n\omega(c_i)\le 0.$$
\end{enumerate}
\defende

The relative coloring test has implications regarding asphericity, but they are not as immediate as in the classical coloring test setting. Given a 2-complex pair $(L,K)$. A spherical diagram $f\colon S\to L$ is called {\em thin relative to $K$} if for every vertex $v\in S$ there exists a 2-cell $d\in S$ with $v$ in its boundary, so that $f(d)\in L-K$.

\satzanf\label{thm:relcolasph}
Let $(L,K)$ be a 2-complex pair with a zero/one angle structure that satisfies the relative coloring test and $\kappa(d)\le 0$ for all 2-cells $d\in L$. If $f\colon S\to L$ is a thin spherical diagram, then it is not vertex reduced.
\satzende

\bewanf Assume $f\colon S\to L$ is a vertex reduced thin spherical diagram. Pull the zero/one-angle structure back to $S$. Then we have $\kappa(d)\le 0$ for every 2-cell $d$ of $S$. Let $v$ be a vertex in $S$. Let $c_1\cdots c_n$ be the corners that make up the link of $v$. Then $f(c_1)\cdots f(c_n)$ is a homology reduced cycle in $lk(f(v),L)$ not entirely contained in $lk(f(v),K)$. It follows that $\kappa(v)\le 0$. We obtain 
$$4=2\chi(S)=\sum_{v\in S}\kappa(v)+\sum_{d\in S}\kappa(d)\le 0$$
and have reached a contradiction. \qed

Under certain conditions this suffices to show that $(L,K)$ is relatively VA. See \cite{HR21}. The relative coloring test also has implications as far as non-positive immersions is concerned. This is work in progress.


\section{The local bi-forest property}

If $\Lambda $ is a graph and $x_1,\ldots ,x_n\in\Lambda$ are vertices then we denote by $\Lambda (x_1,\ldots ,x_n)$ the full subgraph on the vertices $x_1,\ldots ,x_n$.
Let $K$ be a standard 2-complex (i.e. a combinatorial 2-complex with only one vertex), and let $lk(K)$ be the link at its vertex. The link is an undirected graph, we refer to its edges as {\em corners}, because they arise from the corners in 2-cells. If $x_1,\ldots ,x_n$ are the edges of $K$, then the vertices of $lk(K)$ are $\{ x_1^{+},\ldots ,x_n^{+},x_1^{-},\ldots ,x_n^{-}\}$, where $x_i^+$ is located near the start and $x_i^-$ is located near the end of the edge $x_i$, $1\le i\le n$.  We write $x_i^{\epsilon_i}x_j^{\epsilon_j}$, $\epsilon_i, \epsilon_j\in \{ +,-\}$, for the corner that connects $x_i^{\epsilon_i}$ and $x_j^{\epsilon_j}$. This is somewhat sloppy because $lk(K)$ is not a simplicial graph. We will be careful and avoid confusion by providing context whenever necessary. Let $\Lambda = lk(K)$, $\Lambda^+=\Lambda(x_1^{+} ,\ldots ,x_n^{+})$, and $\Lambda^-=\Lambda(x_1^{-} ,\ldots ,x_n^{-})$. If $\epsilon\in \{+, -\}$, then

\begin{equation*}
-\epsilon=
    \begin{cases}
        + & \text{if } \epsilon=-\\
        - & \text{if } \epsilon=+\\
    \end{cases}
\end{equation*}

\defanf\label{dLF} Let $K$ be a standard 2-complex with edges $x_1,\ldots ,x_n$, and let $\Lambda =lk(K)$.

$K$ has the {\em lbf-property (local bi-forest-property)} if
there is a choice of $\epsilon_i\in\{ \pm\}$ for all $1\le i\le n$ such that
$\Lambda(x_1^{\epsilon_1} ,\ldots ,x_n^{\epsilon_n})$ and $\Lambda(x_1^{-\epsilon_1} ,\ldots ,x_n^{-\epsilon_n})$ are forests.

$K$ has the {\em strong lbf-property} if $\Lambda^+$ and $\Lambda^-$ are forests.
\defende

Observe that for every choice of the $\epsilon_i$ the subgraphs $\Lambda(x_1^{\epsilon_1} ,\ldots ,x_n^{\epsilon_n})$ and $\Lambda(x_1^{-\epsilon_1} ,\ldots ,x_n^{-\epsilon_n})$ are disjoint and every vertex of $\Lambda$ belongs to exactly one of the two.

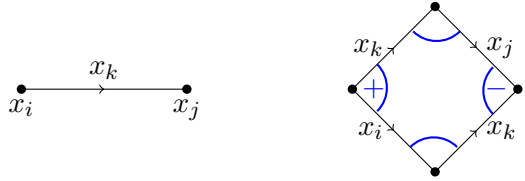
\begin{figure}[ht]\centering
\begin{tikzpicture}[scale=1.1]
\fill (0,1) circle (1.6pt); \node[below] at (0,1) {$x_i$};
\fill (2,1) circle (1.6pt); \node[below] at (2,1) {$x_j$};
\begin{scope}[decoration={markings, mark=at position 0.5 with {\arrow{>}}}]
\draw [postaction={decorate}] (0,1) -- (2,1) node[midway, above]{$x_k$};
\end{scope}

\fill (4,1) circle (1.6pt); \node[right,blue] at (4,1) {$+$};
\fill (5,0) circle (1.6pt);
\fill (5,2) circle (1.6pt);
\fill (6,1) circle (1.6pt); \node[left,blue] at (6,1) {$-$};
\begin{scope}[decoration={markings, mark=at position 0.5 with {\arrow{>}}}]
\draw [postaction={decorate}] (4,1) -- (5,0) node[midway, left]{$x_i$};
\draw [postaction={decorate}] (4,1) -- (5,2) node[midway, left]{$x_k$};
\draw [postaction={decorate}] (5,0) -- (6,1) node[midway, right]{$x_k$};
\draw [postaction={decorate}] (5,2) -- (6,1) node[midway, right]{$x_j$};
\end{scope}

\draw [blue,line width=0.8pt] (4.3,1.3) arc (45:-45:0.4);
\draw [blue,line width=0.8pt] (4.7,0.3) arc (135:45:0.4);
\draw [blue,line width=0.8pt] (5.3,1.7) arc (315:225:0.4);
\draw [blue,line width=0.8pt] (5.7,0.7) arc (225:135:0.4);

\end{tikzpicture}
   \caption{An edge $e$ of $\Gamma$ contributes four corners to $lk(K(\Gamma))$: A positive corner $c_e^+=x_i^+x_k^+$, a negative corner $c_e^-=x_k^-x_j^-$, and two mixed corners, $x_i^-x_k^+$ and $x_k^-x_j^+$.}
   \label{fig:4corners}
\end{figure}

Let $\Gamma$ be a LOG and $e$ an edge of $\Gamma$. Let $s(e)=x_i$, $t(e)=x_j$, and $\lambda(e)=x_k$. Then the 2-cell $d_e$ of $K(\Gamma)$ is attached along the path $x_ix_kx_j^{-1}x_k^{-1}$ and contributes four corners to $lk(K(\Gamma))$. See Figure \ref{fig:4corners}.

\satzanf\label{slfwt} A LOF-complex $K(\Gamma)$ that has the lbf-property admits a zero/ one-angle structure so that the coloring test is satisfied. It follows that $K(\Gamma)$ is DR and therefore aspherical. Furthermore, $K(\Gamma)$ has non-positive immersions which implies that $G(\Gamma)$ is locally indicable. 
\satzende

\bewanf Let $x_1,\ldots ,x_n$ be the edges of $K=K(\Gamma)$ and let $\Lambda=lk(K)$. Assume that $\Lambda_1=\Lambda(x_1^{\epsilon_1} ,\ldots ,x_n^{\epsilon_n})$ and $\Lambda_2=\Lambda(x_1^{-\epsilon_1} ,\ldots ,x_n^{-\epsilon_n})$ are forests.
Assign to every corner of $\Lambda_1\cup\Lambda_2$ angle 0 and all other corners in $\Lambda $ angle 1. Every edge of weight 1 connects a vertex of $\Lambda_1$ to a vertex of $\Lambda_2$. Since $\Lambda_1$ and $\Lambda_2$ are disjoint this implies the cycle condition 2 of the coloring test.

It remains to show that the curvature condition 1 of the coloring test holds. We need to show that every 2-cell of $K$ has 2 corners with angle 0. Assume a 2-cell $d\in K$ is attached along the path $x_ix_kx_j^{-1}x_k^{-1}$. There are 4 cases to consider:
\ben
\item $\epsilon_i=\epsilon_k=\epsilon_j$: The corners $x_i^{+ }x_k^{+ }$ and $x_k^{- }x_j^{- }$ of $d$ have angle 0.
\item $\epsilon_i=-\epsilon_k=\epsilon_j$: The corners $x_i^{- }x_k^{+ }$ and $x_k^{- }x_j^{+ }$ of $d$ have angle 0.
\item $\epsilon_i=\epsilon_k=-\epsilon_j$: The corners $x_i^{+ }x_k^{+ }$ and $x_k^{- }x_j^{+ }$ of $d$ have angle 0.
\item $\epsilon_i=-\epsilon_k=-\epsilon_j$: The corners $x_i^{- }x_k^{+ }$ and $x_k^{- }x_j^{- }$ of $d$ have angle 0.
\een
This show that the coloring test holds. All other statements in the Theorem follow from Theorem \ref{thm:wise}.\qed

Here is one of the main results of this paper. 

\satzanf\label{thm:main} Let $\Gamma$ be a reduced injective LOF such that all sub-LOTs are boundary reduced.  Then the LOF-complex $K(\Gamma)$ has the lbf-property. \satzende

We give a proof in a later section.

\section{Reorientations}

\defanf Given a LOG $\Gamma$. A LOG $\Gamma_{\rho}$ is a {\em reorientation} of $\Gamma$ if it is obtained from $\Gamma$ by reversing the direction of some of its edges. 

A reorientation $\Gamma_{\rho}$ is a {\em block reorientation} if, whenever one edge with label $x$ is reversed, all edges with label $x$ are reversed.
\defende

\lemanf\label{lreo1} Let $\Gamma$ be a LOF with vertex set $\{ x_1,\ldots ,x_n\}$ and let $\Gamma_{\rho}$ be the block reorientation obtained by reversing directions of all edges labeled $x_j$. Then the transposition $x_j^+\leftrightarrow x_j^-$ defines a graph isomorphism 
$$\phi_j\colon \Lambda=lk(K(\Gamma))\to \Lambda_{\rho}=lk(K(\Gamma_{\rho})).$$ In particular
$\Lambda(x_1^{\epsilon_1} ,\ldots, x_j^{\epsilon_j},\ldots, x_n^{\epsilon_n})$ is isomorphic to $\Lambda_{\rho}(x_1^{\epsilon_1} ,\ldots, x_j^{-\epsilon_j},\ldots, x_n^{\epsilon_n})$. 
\lemende

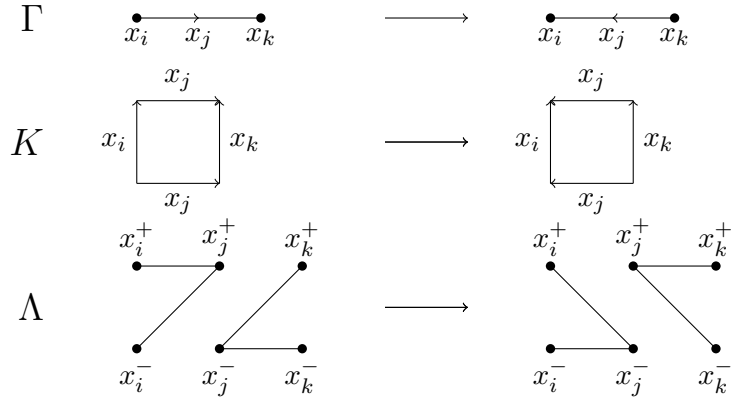
\begin{figure}[ht]\centering
\begin{tikzpicture}[scale=1.1]
\fill (0,0) circle (1.6pt); \node[below] at (0,0) {$x_i^-$};
\fill (1,0) circle (1.6pt); \node[below] at (1,0) {$x_j^-$};
\fill (2,0) circle (1.6pt); \node[below] at (2,0) {$x_k^-$};
\fill (0,1) circle (1.6pt); \node[above] at (0,1) {$x_i^+$};
\fill (1,1) circle (1.6pt); \node[above] at (1,1) {$x_j^+$};
\fill (2,1) circle (1.6pt); \node[above] at (2,1) {$x_k^+$};

 \node[left] at (-1,0.5) {\Large $\Lambda$};
 \node[left] at (-1,2.5) {\Large $K$};
 \node[left] at (-1,4) {\Large $\Gamma$};

\draw (0,0) -- (1,1) -- (0,1);
\draw (2,0) -- (1,0) -- (2,1);

\begin{scope}[decoration={markings, mark=at position 1 with {\arrow{>}}}]
\draw [postaction={decorate}] (0,2) -- (1,2) node[midway, below]{$x_j$};
\draw [postaction={decorate}] (1,2) -- (1,3) node[midway, right]{$x_k$};
\draw [postaction={decorate}] (0,3) -- (1,3) node[midway, above]{$x_j$};
\draw [postaction={decorate}] (0,2) -- (0,3) node[midway, left]{$x_i$};
\draw [postaction={decorate}] (3,2.5) -- (4,2.5);
\draw [postaction={decorate}] (3,0.5) -- (4,0.5);
\draw [postaction={decorate}] (3,4) -- (4,4);
\end{scope}

\fill (5,0) circle (1.6pt); \node[below] at (5,0) {$x_i^-$};
\fill (6,0) circle (1.6pt); \node[below] at (6,0) {$x_j^-$};
\fill (7,0) circle (1.6pt); \node[below] at (7,0) {$x_k^-$};
\fill (5,1) circle (1.6pt); \node[above] at (5,1) {$x_i^+$};
\fill (6,1) circle (1.6pt); \node[above] at (6,1) {$x_j^+$};
\fill (7,1) circle (1.6pt); \node[above] at (7,1) {$x_k^+$};

\draw (5,0) -- (6,0) -- (5,1);
\draw (7,0) -- (6,1) -- (7,1);

\begin{scope}[decoration={markings, mark=at position 1 with {\arrow{>}}}]
\draw [postaction={decorate}] (6,2) -- (5,2) node[midway, below]{$x_j$};
\draw [postaction={decorate}] (6,2) -- (6,3) node[midway, right]{$x_k$};
\draw [postaction={decorate}] (6,3) -- (5,3) node[midway, above]{$x_j$};
\draw [postaction={decorate}] (5,2) -- (5,3) node[midway, left]{$x_i$};
\draw [postaction={decorate}] (3,2.5) -- (4,2.5);
\draw [postaction={decorate}] (3,0.5) -- (4,0.5);
\end{scope}


\begin{scope}[decoration={markings, mark=at position 0.5 with {\arrow{>}}}]
\draw [postaction={decorate}] (0,4) -- (1.5,4) node[midway, below]{$x_j$};
\end{scope}
\fill (0,4) circle (1.6pt); \node[below] at (0,4) {$x_i$};
\fill (1.5,4) circle (1.6pt); \node[below] at (1.5,4) {$x_k$};
\begin{scope}[decoration={markings, mark=at position 0.5 with {\arrow{>}}}]
\draw [postaction={decorate}] (6.5,4) -- (5,4) node[midway, below]{$x_j$};
\end{scope}
\fill (5,4) circle (1.6pt); \node[below] at (5,4) {$x_i$};
\fill (6.5,4) circle (1.6pt); \node[below] at (6.5,4) {$x_k$};
\end{tikzpicture}
\caption{\label{areorient} Reorienting an edge.}
\end{figure}

\bewanf The transposition of the vertices $x_j^+$ and $x_j^-$ is a bijection from the vertex set of $\Lambda$ to the vertex set of $\Lambda_{\rho}$, which are the same sets. We have to check that if $x_p^{\alpha_p}x_q^{\alpha_q}$, $\alpha_p, \alpha_q\in \{ +,-\}$, is a corner of $\Lambda$, then $\phi_j(x_p^{\alpha_p})\phi_j(x_q^{\alpha_q})$ is a corner of $\Lambda_{\rho}$. This is straightforward and one can do this case by case:
\begin{enumerate}
\item $c$ is a corner in a 2-cell $d_e$, where none of $s(e), t(e), \lambda(e)$ is equal to $x_j$.
\item $c$ is a corner in a 2-cell $d_e$, where $s(e)$ or $t(e)$ is $x_j$, but not $\lambda(e)$.
\item $c$ is a corner in a 2-cell $d_e$, such that $\lambda(e)=x_j$; note that $e$ is an edge that gets reversed.
\end{enumerate}
In case 1 the 2-cell $d_e$ exists also in $K(\Gamma_{\rho})$ and the corners are mapped identically. In case 2 $d_e$ exists also in $K(\Gamma_{\rho})$ and the map $\phi_j$ interchanges the two corners that involve $x_j^{\pm}$. We leave the details to the reader, but will do case 3 carefully. Let $e$ be an edge of $\Gamma$ such that $s(e)=x_i$, $t(e)=x_k$, and $\lambda(e)=x_j$.
Then $d_e$ contributes the corners $x_i^-x_j^+, x_j^-x_k^-, x_k^+x_j^-, x_j^+x_i^+$, and $d_{\bar e}$ contributes the corners $\phi_j(x_i^-x_j^+)=x_i^-x_j^-$, $\phi_j(x_j^-x_k^-)=x_j^+x_k^-$, $\phi_j(x_k^+x_j^-)=x_k^+x_j^+$, $\phi_j(x_j^+x_i^+)=x_j^-x_i^+$.
The situation is depicted in Figure \ref{areorient}. \qed

\lemanf\label{lem:notedge} Let $\Gamma$ be a LOF with vertices $\{ x_1,\ldots, x_n\}$ and assume that $x_1$ does not occur as edge label. Let $\Lambda=lk(K(\Gamma))$. Then the transposition $x_1^+\leftrightarrow x_1^-$ defines a graph isomorphism $\phi\colon \Lambda\to \Lambda$. In particular  $\Lambda(x_1^+,x_2^{\epsilon_2},\ldots, x_n^{\epsilon_n})$ is isomorphic to $\Lambda(x_1^-,x_2^{\epsilon_2},\ldots, x_n^{\epsilon_n})$.
\lemende

\bewanf The transposition of the vertices $x_1^+$ and $x_1^-$ is a permutation of the vertex set of $\Lambda$. We have to check that if $x_p^{\alpha_p}x_q^{\alpha_q}$, $\alpha_p, \alpha_q\in \{ +,-\}$, is a corner of $\Lambda$, then $\phi(x_p^{\alpha_p})\phi(x_q^{\alpha_q})$ is a corner of $\Lambda_{\rho}$. This is straightforward and one can do this case by case:
\begin{enumerate}
\item $c$ is a corner in a 2-cell $d_e$, where none of $s(e), t(e), \lambda(e)$ is equal to $x_1$.
\item $c$ is a corner in a 2-cell $d_e$, where $s(e)$ or $t(e)$ is $x_1$. 
\end{enumerate}
In case 1 the corners are mapped identically. Let us do case 2 in more detail. Let $e$ be an edge of $\Gamma$ such that $s(e)=x_1$, $t(e)=x_k$, and $\lambda(e)=x_j\ne x_1$. Then $d_e$ contributes the corners $x_1^-x_j^+, x_j^-x_k^-, x_k^+x_j^-, x_j^+x_1^+$. We have  $\phi(x_1^-x_j^+)=x_1^+x_j^+$, $\phi(x_j^-x_k^-)=x_j^-x_k^-$, $\phi(x_k^+x_j^-)=x_k^+x_j^-$, $\phi(x_j^+x_1^+)=x_j^+x_1^-$. Thus $\phi$ permutes two corners of $d_e$ and fixes the other two. The case where $t(e)=x_1$ is treated similarly. \qed

\satzanf Let $\Gamma$ be a LOF. 
\begin{enumerate}
    \item The lbf-property is invariant under block reorientation.
    \item If $\Gamma$ is injective, then the lbf-property is invariant under reorientation.
    \end{enumerate}
\satzende

\bewanf This is an immediate consequence of Lemma \ref{lreo1}. \qed
\satzanf\label{sinj} An injective LOF $\Gamma$ has the lbf-property if and only if there is a reorientation $\Gamma_{\rho}$ that has the strong lbf-property.
\satzende

\bewanf In the light of Lemma \ref{lreo1} all we have to show is that if $\Gamma$ has the lbf-property then there exists a reorientation that has the strong lbf-property. Let $\{ x_1,\ldots, x_n\}$ be the vertices of $\Gamma$ and assume w.l.o.g.\ $x_1,\ldots, x_k$ do not occur as edge labels, but $x_{k+1},\ldots, x_n$ do. Let $\Lambda=lk(K(\Gamma))$ and assume that $\Lambda(x_1^{\epsilon_1},\ldots x_n^{\epsilon_n})$ and $\Lambda(x_1^{-\epsilon_1},\ldots x_n^{-\epsilon_n})$ are forests. We use Lemma \ref{lreo1} to reorient the edges $e$ of $\Gamma$ for which $\lambda(e)=x_i$ and $\epsilon_i=-$ to produce a reorientation $\Gamma_{\rho}$ so that 
$\Lambda_{\rho}(x_1^{\epsilon_1},\ldots, x_k^{\epsilon_k}, x_{k+1}^+,\ldots x_n^+)$ and $\Lambda_{\rho}(x_1^{\epsilon_1},\ldots, x_k^{\epsilon_k}, x_{k+1}^-,\ldots x_n^-)$ are forests. By Lemma \ref{lem:notedge} $\Lambda_{\rho}(x_1^+,\ldots, x_n^+)$ and 
$\Lambda_{\rho}(x_1^-, \ldots, x_n^-)$ are forests. \qed

\section{The Selection Graph}

For a graph $\Lambda$ we denote its vertex set by $V(\Lambda)$ and its edge set by $E(\Lambda)$. Furthermore, if $E_0$ is a subset of the edge set $E(\Lambda )$, we denote by $\langle E_0 \rangle$ the subgraph of $\Lambda$ spanned by $E$. If $\Lambda$ is a directed graph and $e$ is an edge with $s(e)=x$ and $t(e)=y$, we sometimes write $e=(x\to y)$. 

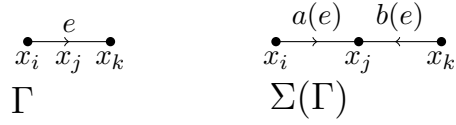
\begin{figure}[ht]\centering
\begin{tikzpicture}[scale=1.1]
\fill (0,0) circle (1.6pt); \node[below] at (0,0) {$x_i$};
\fill (1,0) circle (1.6pt); \node[below] at (1,0) {$x_k$};
\fill (3,0) circle (1.6pt); \node[below] at (3,0) {$x_i$};
\fill (4,0) circle (1.6pt); \node[below] at (4,0) {$x_j$};
\fill (5,0) circle (1.6pt); \node[below] at (5,0) {$x_k$};

 \node[left] at (0.2,-0.7) {\Large $\Gamma$};
 \node[left] at (4,-0.7) {\Large $\Sigma(\Gamma )$};
 \node[above] at (0.5,0) {$e$};
\node[above] at (3.5,0) {$a(e)$};
\node[above] at (4.5,0) {$b(e)$};

\begin{scope}[decoration={markings, mark=at position 0.5 with {\arrow{<}}}]
\draw [postaction={decorate}] (4,0) -- (5,0);
\end{scope}

\begin{scope}[decoration={markings, mark=at position 0.5 with {\arrow{>}}}]
\draw [postaction={decorate}] (0,0) -- (1,0) node[midway, below]{$x_j$};
\draw [postaction={decorate}] (3,0) -- (4,0);
\end{scope}

\end{tikzpicture}
\caption{\label{aKtoS} From $\Gamma$ to $\Sigma(\Gamma )$.}
\end{figure}

Let $\Gamma $ be a LOF and $K(\Gamma)$ be the corresponding 2-complex. We define the {\em selection-graph} $\Sigma(\Gamma )$, a directed graph, as follows: Its vertex set coincides with the vertex set of $\Gamma$ and each edge $e$ of $\Gamma$ gives rise to two edges $a(e)=(s(e)\to \lambda(e))$ and $b(e)=(t(e)\to \lambda(e))$ in $\Sigma(\Gamma)$. (see Figure \ref{aKtoS}).
If $\Gamma$ is compressed, then every edge in $\Sigma(\Gamma )$ has two distinct vertices. $\Sigma(\Gamma )$ may contain multiple edges between two vertices.\\

We immediately make the observation that reorientation of $\Gamma$ has no effect on the selection graph: If $\Gamma_{\rho}$ is a reorientation of $\Gamma$ then $\Sigma(\Gamma)=\Sigma(\Gamma_{\rho})$; however, if en edge $e$ is reversed to $\bar e$, then $a(\bar e)=b(e)$ and $b(\bar e)=a(e)$. 
 
The selection graph $\Sigma(\Gamma)$ is related to both $\Gamma$ and $lk(K(\Gamma))$. Its main point, which we will explain in detail below, is the following: Given an admissible partition into (black and white) subgraphs $\Sigma(\Gamma)=U_b\cup U_w$ it \lq\lq selects\rq\rq\  an orientation $\Gamma_{\rho}$ of $\Gamma$ so that $lk^+(K(\Gamma_{\rho}))=U_b$ and $lk^-(K(\Gamma_{\rho}))=U_w$. Thus, if $\Sigma(\Gamma)$ can be admissibly partitioned into forests, $lk(K(\Gamma_{\rho}))$ and therefore $lk(K(\Gamma))$ is a bi-forest. 
 
Let $\Gamma'$ be the barycentric subdivision of $\Gamma$. If $e$ is an edge in $\Gamma$ let $m_e$ be its midpoint. We orient the edges in $\Gamma'$ to always point towards the $m_e$'s. For $e$ an edge of $\Gamma$ denote by $e^+=(s(e)\to m_e)$ and $e^-=(t(e)\to m_e)$. Note that we have a continuous map 
$$\gamma\colon \Gamma'\to \Sigma(\Gamma)$$
defined by $\gamma(x)=x$, if $x$ is a vertex of $\Gamma$, and $\gamma(m_e)=\lambda(e)$, $\gamma(e^+)=a(e), \gamma(e^-)=b(e)$. $\gamma$ is a bijection on edge sets. Let $\Gamma'_+$ be the subgraph of $\Gamma'$ on the edge set $e^+$, $e\in E(\Gamma)$. Similarly define $\Gamma'_-$. We have maps
$$\alpha_+\colon \Gamma'_+ \to lk^+(K(\Gamma))\  \mbox{and}\ \alpha_-\colon \Gamma'_- \to lk^-(K(\Gamma))$$
defined by $\alpha_+(e^+)=s(e)^+\lambda(e)^+$ and $\alpha_-(e^-)=t(e)^-\lambda(e)^-$. Note that both $\alpha_+$ and $\alpha_-$ are bijections on the edge sets. 

Finally, note that $\Sigma(\Gamma)$ is obtained from $lk^+(K(\Gamma))\cup lk^-(K(\Gamma))$ by identifying $x^+$ and $x^-$ to $x$, for every vertex $x$ of $\Gamma$. Thus we have a quotient map
$$\beta\colon lk^+(K(\Gamma))\cup lk^-(K(\Gamma))\to \Sigma(\Gamma).$$
In detail $\beta$ is described in the following way: $\beta(x^+)=\beta(x^-)=x$; If $e$ is an edge in $\Gamma$ then $\beta(c^+_e)=\beta(s(e)^+\lambda(e)^+)=(s(e)\to \lambda(e))$, and $\beta(c^-_e)=\beta(t(e)^-\lambda(e)^-)=(t(e)\to \lambda(e))$.
Let $\Sigma^+(\Gamma)=\beta(lk^+(\Gamma))$. $\beta$ is a bijection on edges and defines graph isomorphisms $lk^+(\Gamma)\to \Sigma^+(\Gamma)$ and $lk^-(\Gamma)\to \Sigma^-(\Gamma)$. We also note that $\gamma$ factors in the following ways
$$\gamma=\beta\circ \alpha_+ \colon \Gamma'_+\to lk^+(K(\Gamma))\to \Sigma^+(\Gamma)$$ and
$$\gamma=\beta\circ \alpha_- \colon \Gamma'_-\to lk^-(K(\Gamma))\to \Sigma^-(\Gamma).$$
A partition of $E(\Gamma')=E_b\cup E_w$, black and white edges, is {\em admissible} if edges $e^+$  and  $e^-$ have different color for every edge $e\in\Gamma$. A partition $E(\Sigma)=E_b\cup E_w$ is {\em admissible } if edges $a(e)$ and $b(e)$ have different color for each edge $e$. Note that the map $\gamma$ provides a bijection between admissible partitions of $\Gamma'$ and $\Sigma$.

\lemanf\label{ladre} Let $\Gamma$ be a LOF with selection graph $\Sigma (\Gamma)$ and suppose $E(\Sigma)=E_b\cup E_w$ is an admissible partition. Then there exists a reorientation $\Gamma_{\rho}$ of $\Gamma$ so that $E_b=E(\Sigma^+(\Gamma_{\rho}))$ and
$E_w=E(\Sigma^-(\Gamma_{\rho}))$. In particular
$$\beta\colon lk^+(K(\Gamma_{\rho}))\to \Sigma^+(\Gamma_{\rho})=\langle E_b\rangle,\ \beta\colon lk^-(K(\Gamma_{\rho}))\to \Sigma^-(\Gamma_{\rho})=\langle E_w\rangle$$
are graph isomorphisms.
\lemende

\bewanf
The admissible partition for $\Sigma(\Gamma)$ gives an admissible partition for $E(\Gamma')=\gamma^{-1}(E_b)\cup \gamma^{-1}(E_w)$. Now reorient $\Gamma$ to $\Gamma_{\rho}$ so that ${\Gamma'_{\rho}}_+=\langle \gamma^{-1}(E_b)\rangle$. Then 
${\Gamma'_{\rho}}_-=\langle \gamma^{-1}(E_w)\rangle$. Now
$$\gamma({\Gamma'_{\rho}}_+)=\langle E_b \rangle\  \mbox{and}\  
\gamma({\Gamma'_{\rho}}_-)=\langle E_w \rangle.$$
We have factorizations
$$\gamma=\beta\circ\alpha_+\colon {\Gamma'_{\rho}}_+\to lk^+(K(\Gamma_{\rho}))\to \Sigma^+(\Gamma_{\rho})=\langle E_b\rangle$$ and 
$$\gamma=\beta\circ\alpha_-\colon {\Gamma'_{\rho}}_-\to lk^-(K(\Gamma_{\rho}))\to \Sigma^-(\Gamma_{\rho})=\langle E_w\rangle.$$
Since the rightmost maps, the restrictions of $\beta$, are isomorphisms, we obtain the desired result. \qed

\newpage

Let us give a pictorial proof of the Lemma \ref{ladre} (see Figure \ref{fig:KtoS}).
\begin{enumerate}
    \item Draw the LOF $\Gamma$.
    \item Draw the barycentric subdivision $\Gamma'$, and draw $\Sigma(\Gamma)$ by replacing the midpoints $m_e$ by $\lambda(e)$. $\Sigma(\Gamma)$ ``looks'' like $\Gamma'$, but vertices with the same name have to be identified. Best to not actually do these identifications in your drawing since it would only obscure things.
    \item We have given an admissible partition $E(\Sigma)=R\cup B$, red and blue edges (black and white is hard to draw). Put the colors into your drawing. You now see an admissible partition of $E(\Gamma')$ into red and blue edges (half edges of $\Gamma$). Now reorient $\Gamma$ to $\Gamma_{\rho}$ so that if $\bar e$ is an edge in $\Gamma_{\rho}$ then $\bar e^+\in R$. 
    \item Now
    $${\Gamma'_{\rho}}_+\stackrel{\alpha_+}{\to}lk^+(K(\Gamma_{\rho}))\stackrel{\beta}{\to}\langle R \rangle$$ and the latter map is an isomorphism.
\end{enumerate}
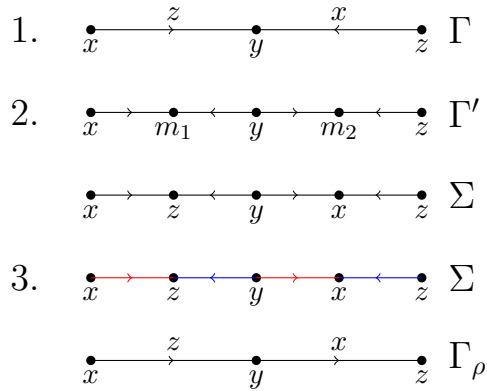
\begin{figure}[ht]\centering
\begin{tikzpicture}[scale=1.1]
\fill (0,0) circle (1.6pt); \node[below] at (0,0) {$x$};
\fill (0,1) circle (1.6pt); \node[below] at (0,1) {$x$};
\fill (0,2) circle (1.6pt); \node[below] at (0,2) {$x$};
\fill (0,3) circle (1.6pt); \node[below] at (0,3) {$x$};
\fill (0,4) circle (1.6pt); \node[below] at (0,4) {$x$};

\fill (2,0) circle (1.6pt); \node[below] at (2,0) {$y$};
\fill (2,1) circle (1.6pt); \node[below] at (2,1) {$y$};
\fill (2,2) circle (1.6pt); \node[below] at (2,2) {$y$};
\fill (2,3) circle (1.6pt); \node[below] at (2,3) {$y$};
\fill (2,4) circle (1.6pt); \node[below] at (2,4) {$y$};

\fill (4,0) circle (1.6pt); \node[below] at (4,0) {$z$};
\fill (4,1) circle (1.6pt); \node[below] at (4,1) {$z$};
\fill (4,2) circle (1.6pt); \node[below] at (4,2) {$z$};
\fill (4,3) circle (1.6pt); \node[below] at (4,3) {$z$};
\fill (4,4) circle (1.6pt); \node[below] at (4,4) {$z$};

\node[right] at (4.2,4) {\Large $\Gamma$};
\node[right] at (4.2,3) {\Large $\Gamma'$};
\node[right] at (4.2,2) {\Large $\Sigma$};
\node[right] at (4.2,1) {\Large $\Sigma$};
\node[right] at (4.2,0) {\Large $\Gamma_\rho$};

\node[left] at (-0.5,4) {\Large $1.$};
\node[left] at (-0.5,3) {\Large $2.$};
\node[left] at (-0.5,1) {\Large $3.$};

\fill (1,1) circle (1.6pt); \node[below] at (1,1) {$z$};
\fill (1,2) circle (1.6pt); \node[below] at (1,2) {$z$};
\fill (1,3) circle (1.6pt); \node[below] at (1,3) {$m_1$};
\fill (3,1) circle (1.6pt); \node[below] at (3,1) {$x$};
\fill (3,2) circle (1.6pt); \node[below] at (3,2) {$x$};
\fill (3,3) circle (1.6pt); \node[below] at (3,3) {$m_2$};

\begin{scope}[decoration={markings, mark=at position 0.5 with {\arrow{<}}}]
\draw [postaction={decorate}] (2,4) -- (4,4) node[midway, above]{$x$};
\draw [postaction={decorate},blue] (1,1) -- (2,1);
\draw [postaction={decorate},blue] (3,1) -- (4,1);
\draw [postaction={decorate}] (1,2) -- (2,2);
\draw [postaction={decorate}] (3,2) -- (4,2);
\draw [postaction={decorate}] (1,3) -- (2,3);
\draw [postaction={decorate}] (3,3) -- (4,3);
\end{scope}

\begin{scope}[decoration={markings, mark=at position 0.5 with {\arrow{>}}}]
\draw [postaction={decorate}] (0,0) -- (2,0) node[midway, above]{$z$};
\draw [postaction={decorate}] (2,0) -- (4,0) node[midway, above]{$x$};
\draw [postaction={decorate}] (0,4) -- (2,4) node[midway, above]{$z$};
\draw [postaction={decorate},red] (0,1) -- (1,1);
\draw [postaction={decorate},red] (2,1) -- (3,1);
\draw [postaction={decorate}] (0,2) -- (1,2);
\draw [postaction={decorate}] (2,2) -- (3,2);
\draw [postaction={decorate}] (0,3) -- (1,3);
\draw [postaction={decorate}] (2,3) -- (3,3);
\end{scope}

\end{tikzpicture}
 \caption{A pictorial proof of Lemma \ref{ladre} }
   \label{fig:KtoS}
\end{figure}

\section{More graph theory: Edmonds' Theorem}

Let $G$ be a directed graph, and let $y$ be a specified vertex in $G$. A {\em branching $B$ in $G$ rooted at $y$} is a spanning tree of $G$ such that for every vertex $v\ne y$, there is exactly one edge in $B$ with $v$ in its boundary which is directed towards $v$.

For a proper subset $X\subseteq V(G)$ we denote by $\delta(X)$ the number of edges that start at a vertex not in $X$ and terminate at a vertex in $X$.

\satzanf\label{thm:edmonds}(Edmonds \cite{Edm73}) There exist $n$ mutually disjoint branchings in $G$ rooted at $y$ if and only if for every subset $X\subseteq V(G)$, $y\notin X$, $\delta(X)\ge n$.
\satzende

A more general result was later shown by Mao-cheng \cite{Mc83}. 
We will apply Edmonds' theorem to the selection graph. For a LOF $\Gamma$ the graph $\Sigma(\Gamma)$ might not be connected. That is why we state the following theorem for LOTs only. 

\satzanf\label{thm:branching} Let $\Gamma$ be a reduced injective LOT such that all sub-LOTs are boundary reduced. Let $y$ be the vertex in $\Gamma$ that does not occur as an edge label. Then there exist two mutually disjoint branchings in $\Sigma(\Gamma)$ rooted at $y$.
\satzende

\lemanf\label{lem:subgraph} Let $\Gamma$ be a reduced injective LOT such that all sub-LOTs are boundary reduced. If $U$ is a subgraph of $\Sigma(\Gamma)$ then $|E(U)|<2|V(U)|-1$.
\lemende

\bewanf Let $T'=\gamma^{-1}(U)$, where $\gamma\colon \Gamma'\to \Sigma(\Gamma)$. Then $T'$ is a subgraph of $\Gamma'$ and hence is a forest. Suppose $T'$ has $l$ components. Then $\chi(T')=l$. Let $y$ be the vertex in $\Gamma$ that does not occur as an edge label. Then $\gamma^{-1}(y)=\{ y \}$. Note that if $x$ is a vertex in $U$ different from $y$ then $\gamma^{-1}(x)=\{ x, m_e\}$, where $m_e$ is the midpoint of the edge $e$ with label $\lambda(e)=x$. Also $|E(T')|=|E(U)|$, because $\gamma$ is a bijection of edges. Thus we have 
$$l=\chi(T')\le 2|V(U)|-|E(U)|,$$
and equality holds if and only if $\gamma^{-1}(x)=\{ x, m_e \}$ for all vertices $x$ of $U$.
We have 
$$|E(U)|\le 2|V(U)|-l.$$
If $l>1$ we are done. So suppose $l=1$, in which case $T'$ is a tree. Assume $1=\chi(T')=2|V(U)|-|E(U)|$. Note that because $\gamma^{-1}(x)=\{ x, m_e \}$, for all $x\in U$, $T'$ contains with every vertex $x$ a subdivided edge of $\Gamma$ with label x, or half an edge of $\Gamma$ with label x. Thus, including missing halfedges (should there be such) and rubbing out edge midpoints produces a sub-LOT $\hat T$ of $\Gamma$. Consider $\hat T'$.
$T'=\hat T'$ is not possible: $\hat T$ would be a sub-LOT where all vertices occur as edge labels, which is impossible because a tree has more vertices than edges.  Therefore $T'<\hat T'$. W.l.o.g. we assume that $T'$ contains the halfedge $e^+$ but not $e^-$ of an edge $e\in\Gamma$. Then, since $T'$ is connected, $m_e$ is a vertex of valency 1 in $T'$. $T'$ contains the halfedge $e^+$ but not $e^-$. Then $t(e)$ is a vertex of $\hat T$ that does not occur as an edge label, so $\hat T$ is not boundary reduced. This contradicts our assumption.  \qed

Proof of Theorem \ref{thm:branching}: Note that since every vertex $x\ne y$ does occur as edge label in $\Gamma$, there are exactly two edges of $\Sigma(\Gamma)$ terminating at $x\ne y$.
Let $X\subset V(\Sigma(\Gamma))$ which does not contain $y$. Let $U$ be the subgraph spanned by $X$. Let $|X|=k$.
Since $X$ does not contain $y$, there are exactly $2k$ edges terminating in $X$. The graph $U$ contains exactly the edges terminating and starting in $X$. Thus $U$ contains $2k-\delta(X)$ edges, since we have to subtract the edges terminating in $X$ but not starting in $X$. Lemma \ref{lem:subgraph} tells us that $$2k-\delta(X)<2k-1$$
and therefore $\delta(X)\ge 2$. The result follows from Theorem \ref{thm:edmonds} since we have shown that for every subset $X\subseteq V(\Sigma(\Gamma))$, $y\notin X$, $\delta(X)\ge 2$. \qed

Let $\Gamma$ be a reduced injective LOT that does contain a sub-LOT $\Gamma_0 $ that is not boundary reduced. Let $X$ be the vertex set of $\Gamma$ and $X_0\subseteq X$ be the vertex set of $\Gamma_0$. Let $x_0$ be the boundary vertex of $\Gamma_0$ that does not occur as an edge label in $\Gamma_0$. Then $\delta(X_0-\{x_0\})=1$. One can show using Edmonds' Theorem \ref{thm:edmonds} that $\Sigma(\Gamma)$ does admit a single branching, but this is not good enough for our purpose.

\section{Proof of Theorem \ref{thm:main}}

We first show the theorem in case $\Gamma$ is a LOT. Let $\Gamma$ be a reduced injective LOT such that all sub-LOTs are boundary reduced, and let $x$ be the vertex in $\Gamma$ that does not occur as an edge label. By Theorem \ref{thm:branching} there are two mutually disjoint branchings $B_1$ and $B_2$ in $\Sigma=\Sigma(\Gamma)$ rooted at $x$. Let $E_i$ be the edge set of $B_i$, $i=1,2$. Note that $E(\Sigma)=E_1\cup E_2$ is an admissible partition. By Lemma \ref{ladre} there exists a reorientation $\Gamma_{\rho}$ so that $lk^+(K_{\rho})=B_1$ and $lk^-(K_{\rho})=B_2$, and hence both are trees. So $K(\Gamma_{\rho})$ has the strong lbf-property. It now follows from Theorem \ref{sinj} that $K$ has the lbf-property.\\

Assume next that $\Gamma$ is an injective reduced LOF with more than one component without boundary reducible sub-LOTs. Assume first that $\Gamma=\Gamma_1\cup \Gamma_2$, where each $\Gamma_i$ is a LOF itself. In that case $K(\Gamma)=K(\Gamma_1)\vee K(\Gamma_2)$. By induction on the number of vertices we can assume that $K(\Gamma_i)$ has the lbf-property, and hence $K(\Gamma)$ has the lbf-property as well.\\

Next assume that $\Gamma$ does not split as above. Then at least one component of $\Gamma$ contains a vertex that is an edge label in a different component. We will see that in this case $\Gamma$ embeds in a reduced injective LOT $\hat\Gamma$ without boundary reducible sub-LOT's. Suppose first $\Gamma=C_1\cup C_2$, where the $C_i$ are the connected components. Let $x_1$ be a vertex of $C_1$ that occurs as an edge label in $C_2$. Let $x_2$ be an arbitrary vertex in $C_2$. Connect $x_1$ to $x_2$ by an edge and label it with a vertex $y\ne x_2$ that does not appear as an edge label in $\Gamma$. This is possible since $\Gamma$ contains exactly 2 vertices which do not appear as edge labels. We obtain a LOT $\hat\Gamma$ with the desired properties. If $\Gamma=C_1\cup\ldots\cup C_k$ we argue by induction to produce $\hat\Gamma$. We have already shown that $K(\hat\Gamma)$ has the lbf-property. Since this property is hereditary and $K(\Gamma)$ is a subcomplex of $K(\hat\Gamma)$, it follows that $K(\Gamma)$ has the lbf-property. \qed

\section{The Relative Case}

Let $K$ be a standard 2-complex and $K_1,\ldots, K_m$ be mutually edge disjoint subcomplexes. Let $X=\{ x_1, \ldots, x_n \}$ be the edge set of $K$ and let $X_i\subseteq X$ be the edge set of $K_i$. Let $X^{\pm}=\{ x_1^+, \ldots, x_n^+, x_1^-, \ldots, x_n^- \}$.  Let $\epsilon\colon X\to X^{\pm}$, $\epsilon(x_i)=x_i^{\epsilon_i}$, $\epsilon_i\in \{ +, - \}$. We have  $\epsilon(X)=\{x_1^{\epsilon_1},\ldots, x_n^{\epsilon_n}\}$. Let $\Lambda=lk(K)$.

\defanf\label{drelLF} 
$(K, K_1\vee\ldots \vee K_m)$ has the {\em relative lbf-property} if
there is a choice of $\epsilon$ such that
$\Lambda(\epsilon(X))$  is a forest relative to $\bigcup_{i=1}^m \Lambda(\epsilon(X_i))$ and $\Lambda(-\epsilon(X))$  is a forest relative to $\bigcup_{i=1}^m \Lambda(-\epsilon(X_i))$.
\defende

\defanf (Quotient LOF) Let $\Gamma$ be a LOF and $\Gamma_1,\ldots, \Gamma_m$ be disjoint sub-LOTs. A {\em quotient LOF} $\bar\Gamma$ is obtained in the following way: For every $1\le i\le m$, choose a vertex $y_i$ from $\Gamma_i$; if $e$ is an edge in $\Gamma$ and $\lambda(e)=z_i$ is a vertex in $\Gamma_i$, then relable $e$ with $y_i$ and collapse $\Gamma_i$ to the vertex $y_i$.
\defende

Note that we have a map $q\colon K(\Gamma)\to K(\bar\Gamma)$. Here is how $\bar\Lambda=lk(K(\bar\Gamma))$ is obtained from $\Lambda=lk(K(\Gamma))$. For every $1\le i\le m$, remove all corners of $\Lambda_i=lk(K(\Gamma_i))$ that are not contained in $\Lambda_i^+ \cup \Lambda_i^-$. For every $1\le i\le m$,  identify all of $\Lambda_i^+$ to $y_i^+$, and all of $\Lambda_i^-$ to $y_i^-$. 

\satzanf\label{thm:rellbf}
Let $\Gamma$ be a reduced LOF, $\Gamma_1,\ldots, \Gamma_m$ be disjoint sub-LOTs, and $\bar \Gamma$ be the quotient LOF. If $K(\bar\Gamma)$ has the lbf-property, then the pair $(K(\Gamma), K(\Gamma_1)\vee\ldots\vee K(\Gamma_m))$ has the relative lbf-property.
\satzende

\bewanf Let $\Lambda=lk(K(\Gamma))$, $\Lambda_i=lk(K(\Gamma_i))$, and $\bar\Lambda=lk(K(\bar\Gamma))$. Let $\bar X$, $X$, $X_i$ the vertex sets of $\bar \Gamma$, $\Gamma$ and $\Gamma_i$, respectively. There is an $\bar\epsilon$ so that both $\bar\Lambda(\bar\epsilon(\bar X))$ and $\bar\Lambda(-\bar\epsilon(\bar X))$ are forests. We construct an $\epsilon$ in the following way: If $x\in X$ and $x$ is not contained in any $X_i$, then define $\epsilon(x)=\bar\epsilon(x)$. If $y_i$ is the vertex of $\bar\Gamma_i$ selected in the quotient process, and $z_i\in X_i$, define $\epsilon(z_i)=\bar\epsilon(y_i)$. Note that $\Lambda_i(\epsilon(X_i))=\Lambda_i^+$ if $\epsilon(y_i)=y_i^+$, and  $\Lambda_i(\epsilon(X_i))=\Lambda_i^-$ if $\epsilon(y_i)=y_i^-$. It follows that $\bar\Lambda(\bar\epsilon(\bar X))$ is obtained from $\Lambda(\epsilon(X))$ by collapsing $\Lambda_i(\epsilon(X_i))$ to $\epsilon(y_i)$, for each $1\le i\le m$. Thus we do have a quotient map
$$q\colon\Lambda(\epsilon(X))\to \bar\Lambda(\bar\epsilon(\bar X)).$$
Since $\bar\Lambda(\bar\epsilon(\bar X))$ is a forest, $\Lambda(\epsilon(X))$ is a forest relative to $\bigcup_{i=1}^m\Lambda_i(\epsilon(X_i))$.  Indeed, let $c_1\cdots c_k$ be a cycle of corners in $\Lambda(\epsilon(X))$ that is not entirely contained in  $\bigcup_{i=1}^m\Lambda_i(\epsilon(X_i))$, and let $q(c_1)\cdots q(c_k)$ be the corresponding cycle in $\bar\Lambda(\bar\epsilon(\bar X))$.
Note that $q(c_j)=\epsilon(y_i)$ if $c_j$ is a corner contained in $\Lambda(\epsilon(X_i))$. Thus $q(c_j)$ is not a corner and does not appear in the corner cycle $q(c_1)\cdots q(c_k)$. Slightly abusing notation we write $q(c_j)=\emptyset$. By our assumption not all $q(c_j)=\emptyset$. Since $\bar\Lambda(\bar\epsilon(\bar X))$ is a forest, there has to be a pair $q(c_i),q(c_j)$ so that $q(c_j)=q(c_i)^{opp}$, the corner $q(c_i)$ with opposite orientation, and $q(c_l)=\emptyset$ for $i<l<j$. Thus $c_j=c_i^{opp}$, and our cycle is not homology reduced. The fact that $\Lambda(-\epsilon(X))$ is a forest relative to $\bigcup_{i=1}^m\Lambda_i(-\epsilon(X_i))$ is shown in the same way.\qed

\satzanf\label{thm:relcolLOT}
Let $\Gamma$ be a reduced LOF, $\Gamma_1,\ldots, \Gamma_m$ be disjoint sub-LOTs, and $\bar \Gamma$ be the quotient LOF. If $K(\bar\Gamma)$ has the lbf-property, then $K(\Gamma)$ admits a zero/one-angle structure so that the pair $(K(\Gamma), K(\Gamma_1)\vee\ldots\vee K(\Gamma_m))$ satisfies the relative coloring test and $\kappa(d)\le 0 $ for all 2-cells $d$ of $K(\Gamma)$.
\satzende

\bewanf
Let $\Lambda=lk(K(\Gamma))$, $\Lambda_i=lk(K(\Gamma_i))$, and $\bar\Lambda=lk(K(\bar\Gamma))$. Let $\bar X$ be the vertex set of $\bar\Gamma$. There is an $\bar\epsilon$ so that both $\bar\Lambda(\bar\epsilon(\bar X))$ and $\bar\Lambda(-\bar\epsilon(\bar X))$ are forests. We construct an $\epsilon$ as in the proof of the previous Theorem \ref{thm:rellbf}. Then $\Lambda(\epsilon(X))$ is a forest relative to $\bigcup_{i=1}^m\Lambda_i(\epsilon(X_i))$ and $\Lambda(-\epsilon(X))$ is a forest relative to $\bigcup_{i=1}^m\Lambda_i(-\epsilon(X_i))$. We now proceed in the same way as in the proof of Theorem \ref{slfwt}: Assign to a corner in $\Lambda(\epsilon(X))\cup \Lambda(-\epsilon(X))$ angle $0$, and to all other corners angle 1. Let $c_1\cdots c_k$ be a homology reduced cycle of corners in $\Lambda$ not entirely contained in $\bigcup_{i=1}^m \Lambda_i$. Since both $\Lambda(\epsilon(X))$ and $\Lambda(-\epsilon(X))$ are relative forests, our cycle can not be contained in $\Lambda(\epsilon(X))$ or $\Lambda(-\epsilon(X))$. Since these two are disjoint, it follows that our cycle must contain two corners $c_i$ and $c_j$ not contained in $\Lambda(\epsilon(X))\cup\Lambda(-\epsilon(X))$. Since their angles $\omega(c_i)=\omega(c_j)=1$, the cycle condition 2 in the relative coloring test is satisfied.

We next check the curvature of 2-cells. First consider a 2-cell $d_e$, where $e$ is an edge in $\Gamma$ that is not an edge in any $\Gamma_i$. Let $\bar e$ be the corresponding edge in $\bar\Gamma$. It was shown in Theorem \ref{slfwt} that $d_{\bar e}$ contains two corners from $\Lambda(\bar\epsilon(\bar X))\cup \Lambda(-\bar\epsilon(\bar X))$. Therefore $d_e$ contains two corners from $\Lambda(\epsilon(X))\cup \Lambda(-\epsilon(X))$ and hence two corners of angle 0. It follows that $\kappa(d_e)\le 0$.
Next assume that $e$ is an edge of $\Gamma_i$, for some $i$. We assume w.l.o.g.\ that $\bar\epsilon(y_i)=+$. Then $\Lambda(\epsilon(X_i))=\Lambda_i^+\subseteq \Lambda(\epsilon(X))$ and $\Lambda(-\epsilon(X_i))=\Lambda_i^-\subseteq \Lambda(-\epsilon(X))$. It follows that both corners $c_e^+$ and $c_e^-$ have angles $0$ and therefore $\kappa(d_e)\le 0$. \qed

\satzanf\label{thm:aspherical} If $\Gamma$ is a reduced injective LOF then $K(\Gamma)$ is VA. It follows that every injective LOF is aspherical. 
\satzende

This result has appeared in our earlier work \cite{HR17} and \cite{HR21}. We give a proof sketch. Assume that $\Gamma$ is reduced and injective. Consider the maximal sub-LOTs $\Gamma_1,\ldots , \Gamma_m$. The generic situation is this: The maximal sub-LOTs are disjoint and the quotient $\bar \Gamma$ is reduced injective and without sub-LOTs. Other situations can be dealt with in an ad hoc fashion. It follows from Theorem \ref{thm:main} that $K(\bar \Gamma)$ has the lbf-property and so the pair $(K(\Gamma)$, $K(\Gamma_1)\vee\ldots\vee K(\Gamma_m))$ satisfies the relative coloring test by Theorem \ref{thm:relcolLOT} and $\kappa(d)\le 0$ for all 2-cells of $K(\Gamma)$. Let $L=K(\Gamma)$ and $J=K(\Gamma_1)\vee\ldots\vee K(\Gamma_m)$. It follows from Theorem \ref{thm:relcolasph} that there are no vertex reduced thin spherical diagrams over $(L,J)$. It turns out, using the concept of thinning expansions (see Definition 3.1 in \cite{HR21}), that this is enough to show that $(L,J)$ is relatively VA.
$J=K(\Gamma_1)\vee\ldots\vee K(\Gamma_m)$, and since each $\Gamma_i$ has fewer vertices than $\Gamma$ we can inductively assume that each $K(\Gamma_i)$ is VA. (A little care is called for here, because the LOTs $\Gamma_i$ might not be boundary reduced. But it is not difficult to see that if a LOT-complex is VA after boundary reductions, it is VA to begin with.) Therefore $J$ is VA. Now $(L,J)$ is relatively VA and $J$ is VA, and that implies that $L$ is VA.

Assume that $\Gamma$ is injective, but not reduced. Since reductions do not change the homotopy type of the LOF complex, it follows that $K(\Gamma)$ is aspherical.\qed

\vspace{1cm}
Department of Mathematics,
Boise State University,
Boise, ID 83725-1555,
USA

email: jensharlander@boisestate.edu

\bigskip
P{\"a}dagogische Hochschule Karlsruhe,
Bismarckstr. 10,
76133 Karlsruhe,
Germany
 
email: rosebrock@ph-karlsruhe.de

\end{document}